\documentclass{proc-l}

\usepackage{amssymb}
\usepackage[cp850]{inputenc}

\newtheorem{theorem}{Theorem}
\newtheorem{lemma}[theorem]{Lemma}

\newtheorem{corollary}[theorem]{Corollary}

\theoremstyle{definition}

\theoremstyle{remark}

\def\Ric{\mathop\mathrm{Ric}\nolimits}
\newcommand{\m}{\mbox{$M$}}

\newcommand{\s}{\mbox{$\Sigma$}}
\newcommand{\Si}{\mbox{${\mathbb S}$}}
\newcommand{\R}{\mbox{${\mathbb R}$}}
\newcommand{\Nn}{\mbox{$-\R \times_{e^t} \m^n$}}

\newcommand{\g}[2]{\mbox{$\langle #1 ,#2 \rangle$}}
\newcommand{\fle}{\mbox{$\rightarrow$}}
\newcommand{\rf}[1]{\mbox{(\ref{#1})}}
\newcommand{\rl}[1]{{~\ref{#1}}}
\newcommand{\nablabar}{\mbox{$\overline{\nabla}$}}
\newcommand{\fs}{\mbox{$\mathcal{C}^\infty(\s)$}}

\def\beq{\begin{equation}}
\def\eeq{\end{equation}}

\begin{document}

\title[Spacelike hypersurfaces in the steady state space]
{Spacelike hypersurfaces with constant mean curvature in the steady state space}
\thanks{This research is a result of the activity developed within the framework of the
Programme in Support of Excellence Groups of the Regi\'{o}n de Murcia, Spain, by Fundaci\'{o}n S\'{e}neca,
Regional Agency for Science and Technology (Regional Plan for Science and Technology 2007-2010).}
\thanks{
Research partially supported by MEC project MTM2007-64504, and Fundaci\'{o}n S\'{e}neca project 04540/GERM/06,
Spain.}
\thanks{
A.L. Albujer was supported by FPU Grant AP2004-4087 from Secretar\'\i a de Estado de Universidades e
Investigaci\'{o}n, MEC Spain.}

\author{Alma L. Albujer}
\address{Departamento de Matem\'{a}ticas, Universidad de Murcia, E-30100 Espinardo, Murcia, Spain}
\email{albujer@um.es}

\author{Luis J. Al\'\i as}
\address{Departamento de Matem\'{a}ticas, Universidad de Murcia, E-30100 Espinardo, Murcia, Spain}
\email{ljalias@um.es}

\subjclass[2000]{53C42, 53C50}

\date{May 2007}


\keywords{de Sitter space, steady state space, spacelike hypersurface, mean curvature, parabolicity}

\begin{abstract}
We consider complete spacelike hypersurfaces with constant mean curvature in the open region of de Sitter space
known as the steady state space. We prove that if the hypersurface is bounded away from the infinity of the ambient
space, then the mean curvature must be $H=1$. Moreover, in the 2-dimensional case we obtain that the only complete
spacelike surfaces with constant mean curvature which are bounded away from the infinity are the totally umbilical flat
surfaces. We also derive some other consequences for hypersurfaces which are bounded away from the future infinity. Finally,
using an isometrically equivalent model for the steady state space, we extend our results to a wider family of spacetimes.
\end{abstract}

\maketitle

\section{The steady state space}
Let $\R_1^{n+2}$ be the $(n+2)$-dimensional Lorentz-Minkowski
space, that is, the real vector space $\R^{n+2}$ endowed with the
Lorentzian metric
\[
\g{}{}=-(dx_0)^2+(dx_1)^2+...+(dx_{n+1})^2
\]
where $(x_0,...,x_{n+1})$ are the canonical coordinates of $\R^{n+2}$. As
it is well known, the hyperquadric
\[
\Si_1^{n+1}=\{x \in \R_1^{n+2}: \g{x}{x}=1\}
\]
consisting of all unit spacelike vectors in $\R_1^{n+2}$ endowed with the
induced metric from $\R_1^{n+2}$ is the de Sitter space. The de Sitter
space is a complete simply connected $(n+1)$-dimensional Lorentzian
manifold with constant sectional curvature one. Therefore, $\Si_1^{n+1}$
can be seen, in Lorentzian geometry, as the equivalent to the Euclidean
sphere.

Take now a non-zero null vector  $a \in \R_1^{n+2}$ past-pointing, that is,
$\g{a}{a}=0$ and $\g{a}{e_0}>0$ where $e_0=(1,0,...,0)$ and consider the
following open region of the de Sitter space
\[
\mathcal{H}^{n+1}=\{x \in \Si_1^{n+1}: \g{x}{a}>0\}.
\]
The open region $\mathcal{H}^{n+1}$ forms the spacetime for the steady state model of the universe proposed by
Bondi and Gold \cite{BG} and Hoyle \cite{Ho}, when looking for a model of the universe which looks the same not only
at all points and in all directions (that is, spatially isotropic and homogeneous), but at all times (cf.
\cite[Section 14.8]{We} and \cite[Section 5.2]{HE}). The steady state space is a
non-complete manifold, being only half of the de Sitter space and having as boundary the null hypersurface
\[
L_0=\{x \in \Si_1^{n+1}: \g{x}{a}=0\}.
\]
The null hypersurface $L_0$ represents the past infinity of $\mathcal{H}^{n+1}$, usually denoted by $\mathcal{J}^{-}$, and
the limit boundary
\[
L_\infty=\{x \in \Si_1^{n+1}: \g{x}{a}=\infty \}
\]
represents its future infinity, usually denoted by $\mathcal{J}^{+}$.

\section{Spacelike hypersurfaces in the steady state space}
A smooth immersion $\psi:\Sigma^n \rightarrow \mathcal{H}^{n+1}$ of an $n$-dimensional connected manifold $\Sigma^n$ is
said to be a spacelike hypersurface if $\psi$ induces a Riemannian metric on $\Sigma$, which as usual is also denoted by
$\g{}{}$. In that case, there exists a unique unitary timelike normal field $N$ globally defined on $\Sigma$ which is
future-directed. Throughout this paper we will refer to $N$ as the future-pointing Gauss map of $\Sigma$. The mean
curvature function of a spacelike hypersurface $\Sigma$ is defined as
\[
H=-\frac{1}{n}\mathrm{tr}(A),
\]
where $A$ stands for the shape operator (or second fundamental form) of $\Sigma$ with respect to its future-pointing
Gauss map $N$. The choice of the sign $-$ in our definition of $H$ is motivated by the fact that in that case the
mean curvature vector is given by $\overrightarrow{H}=HN$. Therefore, $H(p)>0$
at a point $p\in\Sigma$ if and only if $\overrightarrow{H}(p)$ is future-directed.

The steady state space has a special property, it admits a foliation by means of totally umbilical
spacelike hypersurfaces
\[
L_\tau=\{x \in \mathcal{H}^{n+1}: \g{x}{a}=\tau \}, \quad \tau >0,
\]
having constant mean curvature $H=1$ with respect to their future-pointing Gauss map
\[
N_\tau (x)=x-\frac{1}{\tau}a.
\]
It is worth pointing out that each $L_\tau$ is isometric to the Euclidean space
$\R^n$. Actually, if we identify $\R^n$ with the orthogonal complement in
$\R^{n+2}_1$ of the Lorentzian plane spanned by $a$ and $e_0$, then it is easy to
see that the map $\phi_\tau:L_\tau\fle\R^n$ given by
\[
\phi_\tau(x)=x-\frac{\tau+\g{a}{e_0}\g{x}{e_0}}{\g{a}{e_0}^2}a-\frac{\tau}{\g{a}{e_0}}e_0
\]
defines an isometry between $L_\tau$ and $\R^n$.

In this context, we will say that a spacelike hypersurface $\Sigma$ in $\mathcal{H}^{n+1}$ is
\textit{bounded away from the future infinity} if there exists $\overline{\tau}>0$ such that
\[
\psi(\Sigma)\subset \{x \in \mathcal{H}^{n+1}: \g{x}{a} \leq \overline{\tau} \},
\]
and we will say that it is \textit{bounded away from the past infinity} if there exists $\underline{\tau}>0$ such that
\[
\psi(\Sigma)\subset \{x \in \mathcal{H}^{n+1}: \g{x}{a} \geq \underline{\tau} \}.
\]
We will say that $\Sigma$ is \textit{bounded away from the infinity} if it is both bounded away from the past and the
future infinity. In other words, $\Sigma$ is bounded away from the infinity if there exist
$0<\underline{\tau}<\overline{\tau}$ such that $\psi(\Sigma)$ is contained in the slab bounded by $L_{\underline{\tau}}$
and $L_{\overline{\tau}}$.

\begin{lemma}
\label{lemma}
Let $\psi:\Sigma^n\fle\mathcal{H}^{n+1}$ be a complete spacelike hypersurface. If $\Sigma$ is bounded away from the future
infinity of $\mathcal{H}^{n+1}$, then $\Sigma$ is diffeomorphic to $\R^{n}$. In particular, there is no compact
(without boundary) spacelike hypersurface in $\mathcal{H}^{n+1}$.
\end{lemma}
\begin{proof}
We will see that $\Sigma$ is diffeomorphic to $\R^n$ by showing that its orthogonal geodesic projection onto $L_1$ is a
diffeomorphism. Observe that the steady state spacetime $\mathcal{H}^{n+1}$ can be globally parametrized by means of the
diffeomorphism $\Phi:\R\times L_1\fle\mathcal{H}^{n+1}$ given by $\Phi(t,q)=\gamma_{q}(t)$, where
\[
\gamma_q(t)=e^{t}\ q-\sinh{t}\ a
\]
denotes the (future pointing) unitary geodesic orthogonal to $L_1$ through the point $q\in L_1$. Then, the orthogonal
geodesic projection of $\mathcal{H}^{n+1}$ onto $L_1$ assigns to each $x\in\mathcal{H}^{n+1}$ the point $\Psi(x)\in L_1$
such that $x=\Phi(t,\Psi(x))$ for a certain $t$, and it is given by
\[
\Psi(x)=\frac{1}{\g{x}{a}}x+\frac{1}{2}\left(1-\frac{1}{\g{x}{a}^2}\right)a, \quad x\in\mathcal{H}^{n+1}.
\]
Let $\psi:\Sigma^n\fle\mathcal{H}^{n+1}$ be a spacelike hypersurface. The orthogonal geodesic projection of $\Sigma$ onto
$L_1$ is the composition map $\Pi=\Psi\circ\psi:\Sigma\fle L_1$ written as
\[
\Pi=\frac{1}{\g{\psi}{a}}\psi+\frac{1}{2}\left(1-\frac{1}{\g{\psi}{a}^2}\right)a.
\]
A straightforward computation shows that
\[
d\Pi_p(v)=\frac{1}{\g{\psi(p)}{a}}d\psi_p(v)
-\frac{\g{d\psi_p(v)}{a}}{\g{\psi(p)}{a}^2}\psi(p)
+\frac{\g{d\psi_p(v)}{a}}{\g{\psi(p)}{a}^3}a
\]
for every $p\in\Sigma$ and every tangent vector $v\in T_p\Sigma$. Therefore,
\begin{eqnarray}
\label{e2}
\g{d\Pi_p(v)}{d\Pi_p(v)}_{\mathrm{o}}& = & \frac{1}{\g{\psi(p)}{a}^2}\g{d\psi_p(v)}{d\psi_p(v)}+
\frac{\g{d\psi_p(v)}{a}^2}{\g{\psi(p)}{a}^4}\\
\nonumber {} & \geq & \frac{1}{\g{\psi(p)}{a}^2}\g{d\psi_p(v)}{d\psi_p(v)},
\end{eqnarray}
that is,
\[
\Pi^*(\g{}{}_{\mathrm{o}})\geq\frac{1}{\g{\psi}{a}^2}\g{}{}
\]
where $\g{}{}_{\mathrm{o}}$ denotes the flat Euclidean metric on $L_1$ and $\g{}{}$ denotes the Riemannian metric on
$\Sigma$. Since we are assuming that $\Sigma$ is bounded away from the future infinity of $\mathcal{H}^{n+1}$,
it follows from here that
\beq
\label{e4}
\Pi^*(\g{}{}_{\mathrm{o}})\geq\frac{1}{\overline{\tau}^2}\g{}{}
\eeq
for certain positive $\overline{\tau}$.

From \rf{e2} we get that $\Pi$ is a local diffeomorphism. Since $\g{}{}$ is a
complete Riemannian metric on $\Sigma$, the same holds for the homothetic metric
$\widetilde{\g{}{}}=(1/\overline{\tau}^{2})\g{}{}$. Then, \rf{e4} means that the
map
\[
\Pi:(\Sigma,\widetilde{\g{}{}})\fle(L_1,\g{}{}_{\mathrm{o}})
\]
increases the distance, and hence $\Pi$ is a covering map by \cite[Chapter VIII, Lemma 8.1]{KN}.
But $L_1$ being simply connected this means that $\Pi$ is in fact a global diffeomorphism between $\Sigma$ and $L_1$.
\end{proof}

\section{Spacelike hypersurfaces with constant mean curvature}
For a spacelike hypersurface $\psi:\Sigma^n \fle \mathcal{H}^{n+1}$, consider
the function $u \in \fs$ defined as $u(p)=\g{\psi(p)}{a}$, $p \in \Sigma$.
Observe that the gradient of $u$ in $\Sigma$ is
\[
\nabla u =a ^\top,
\]
where $a^\top$ denotes the tangential component of $a$ along $\Sigma$, that is,
\beq
\label{adesc}
a=a^\top-\g{N}{a}N+\g{\psi}{a}\psi=\nabla u-\g{N}{a}N+\g{\psi}{a}\psi.
\eeq
Using Gauss and Weingarten
formulae, we easily obtain
\[
\nabla_X \nabla u=-\g{N}{a}AX-u X
\]
for every $X \in T \Sigma$. Therefore, the Laplacian of the function $u$ on $\Sigma$ is given by
\beq
\label{laplau}
\Delta u=n H \g{N}{a}- n u.
\eeq
From \eqref{adesc}, it is also easy to see that
\beq
\label{normnablau}
\| \nabla u \| ^2=\g{N}{a}^2-u^2
\eeq
where $\|\cdot\|$ denotes the norm of a vector field on $\Sigma$. Recall here that $a$ is a past-pointing null vector and
$N$ is future-pointing, so that $\g{N}{a}>0$ on $\Sigma$.

In order to prove our main results, we will make use of the well-known
Omori-Yau maximum principle \cite{Yau}, recalled next.
\begin{lemma}
\label{OmoriYau} Let $M$ be a complete Riemannian manifold whose Ricci
curvature is bounded from below. If $u \in \mathcal{C}^\infty(M)$ is
bounded from above on $M$ then there exists a sequence of points
$\{p_k\}_{k\in\mathbb{N}}$ in $M$ such that
\[
\lim_{k\rightarrow\infty} u(p_k)=\sup_M u, \quad \|\nabla u(p_k)\|<1/k
\;\;\; \mathrm{and} \;\;\; \Delta u(p_k)<1/k.
\]
\end{lemma}

Now we are ready to state the following result.
\begin{theorem}
\label{deSitter}
Let $\psi:\Sigma^n\fle\mathcal{H}^{n+1}$ be a complete spacelike hypersurface with constant mean curvature $H$.
If $\Sigma$ is bounded away from the infinity of $\mathcal{H}^{n+1}$, then $H=1$ necessarily. Moreover, in the
$2$-dimensional case, there exists $\tau$ such that $\Sigma^2=L_{\tau}$.
\end{theorem}
\begin{proof}
The Gauss equation of $\Sigma$ in $\mathcal{H}^{n+1}$ describes the curvature of $\Sigma$, denoted by $R$,
in terms of its shape operator, and it is given by
\[
\g{R(X,Y)X}{Y}=\g{X}{X}\g{Y}{Y}-\g{X}{Y}^2-\g{AX}{X}\g{AY}{Y}+\g{AX}{Y}^2
\]
being $X,Y \in T\Sigma$. Taking traces here, we obtain that for every $X \in T\Sigma$, $\|X\|=1$,
\begin{eqnarray}
\label{ricbound}
\Ric(X,X) & = & n-1+nH\g{AX}{X}+\g{AX}{AX}\\
\nonumber {} & = & n-1+\left\| AX+\frac{nH}{2}X \right \|^2 -\frac{n^2H^2}{4}\geq n-1-\frac{n^2H^2}{4},
\end{eqnarray}
where $\Ric$ stands for the Ricci curvature of $\Sigma$.
Thus, the Ricci curvature
of $\Sigma$ is bounded from below by the constant $n-1-n^2H^2/4$. By assumption the
function $u$ is bounded (from above and from below). Therefore, applying Lemma
\ref{OmoriYau} to the function $u$, we know that there exists a sequence
$\{p_k\}_{k\in\mathbb{N}}$ in $\Sigma$ such that \beq \label{lim1} \lim_{k
\rightarrow \infty} u(p_k)=\sup_\Sigma u<+\infty, \eeq \beq \label{lim2} \| \nabla
u(p_k) \|^2=\g{N(p_k)}{a}^2 - u^2(p_k) < 1/k^2 \eeq and
\[
\Delta u(p_k)=n(H \g{N(p_k)}{a}-u(p_k))<1/k.
\]
From the last equation, taking into account that $\g{N}{a}>0$ on $\Sigma$, we obtain that
\beq
\label{ineq1}
H<\frac{u(p_k)}{\g{N(p_k)}{a}}+\frac{1}{n\g{N(p_k)}{a}k}.
\eeq
On the other hand, taking limits in \eqref{lim2} and using \eqref{lim1} we also get that
\[
\lim_{k \rightarrow \infty} \g{N(p_k)}{a}=
\lim_{k \rightarrow \infty}u(p_k) = \sup_\Sigma u.
\]
Therefore, taking limits in \eqref{ineq1} we conclude that $H \leq 1$.

In an similar way, we may apply Lemma \ref{OmoriYau} to the function
$-u$, obtaining another sequence $\{q_k\}_{k\in\mathbb{N}}$ in $\Sigma$ such that
\[
\lim_{k \rightarrow \infty} u(q_k)=\inf_\Sigma u>0,
\]
\[
\| \nabla u(q_k) \|^2=\g{N(q_k)}{a}^2 - u^2(q_k) < 1/k^2
\]
and
\[
\Delta u(q_k)=n(H \g{N(q_k)}{a}-u(q_k))>-1/k.
\]
Therefore, we get now that
\[
H>\frac{u(q_k)}{\g{N(q_k)}{a}}-\frac{1}{n\g{N(q_k)}{a}k}
\]
and taking limits again we conclude that $H \geq 1$. Summing up $H=1$. Moreover, when $n=2$ we know by a result due to
Akutagawa \cite{Ak} and Ramanathan \cite{Ra} that $\Sigma$ must be a totally umbilical surface of the foliation $L_\tau$.
For the sake of completeness, we give here a different and direct proof of this assertion in our situation.

Hence, assume that $n=2$. Since $H=1$, \eqref{ricbound} implies that the Gaussian
curvature of $\Sigma$ is non-negative. Then $\Sigma$ is a complete Riemannian surface with non-negative Gaussian
curvature, and by a classical result due to Ahlfors \cite{Ah} and Blanc-Fiala-Huber \cite{Hu}, $\Sigma$ is parabolic, in
the sense that any subharmonic function bounded from above on the surface must be constant.
By \eqref{normnablau} and from the fact that $\g{N}{a}>0$ we obtain that
\[
\g{N}{a}-u  \geq 0
\]
on $\Sigma$. Then, \eqref{laplau} implies that $\Delta u \geq 0$. Therefore $u$ is a subharmonic function bounded
from above on the parabolic surface $\Sigma$, and hence it must be
constant, being $\Sigma$ one of the totally umbilical spacelike surfaces $L_\tau$.
\end{proof}
\begin{corollary}
\label{coro1}
The only complete spacelike surfaces with constant mean curvature in $\mathcal{H}^{3}$ which are bounded away from the
infinity are the totally umbilical surfaces $L_\tau$.
\end{corollary}

As a consequence of our Lemma\rl{lemma} and the proof of Theorem\rl{deSitter} we can also state the following result.
\begin{theorem}
\label{deSitterBIS}
Let $\psi:\Sigma^n\fle\mathcal{H}^{n+1}$ be a complete spacelike hypersurface with constant mean curvature $H$.
If $\Sigma$ is bounded away from the future infinity of $\mathcal{H}^{n+1}$ and has future-pointing mean
curvature vector, then $2\sqrt{n-1}/n\leq H\leq 1$. Moreover, in the $2$-dimensional case,
there exists $\tau$ such that $\Sigma^2=L_{\tau}$.
\end{theorem}
\begin{proof}
In this case the function $u$ is bounded only from above. Therefore, by the first part of the proof of
Theorem\rl{deSitter} we obtain that $H\leq 1$. On the other hand, estimate \rf{ricbound} implies that
$H^2\geq 4(n-1)/n^2$. Otherwise by Bonnet-Myers' theorem we would get that $\Sigma$ is compact, which is not possible by
Lemma\rl{lemma}. But the mean curvature vector field being future-pointing yields $H\geq 2\sqrt{n-1}/n$. Summing up,
$2\sqrt{n-1}/n\leq H\leq 1$. When $n=2$ this gives $H=1$ and then $\Sigma$ must be one of the totally umbilical
spacelike surfaces $L_\tau$.
\end{proof}
\begin{corollary}
\label{coro2}
The only complete spacelike surfaces with constant mean curvature in $\mathcal{H}^{3}$ which are bounded away from the
future infinity and have future-pointing mean curvature vector are the totally umbilical surfaces $L_\tau$.
\end{corollary}

It is worth pointing out that Corollary\rl{coro2} is no longer true for surfaces which are bounded away from the past
infinity of $\mathcal{H}^{3}$.
Actually, in \cite[Corollary 12]{Mo} Montiel constructs complete spacelike hypersurfaces in $\mathcal{H}^{n+1}$ with
constant mean curvature $H>1$ which are bounded away from the past infinity. On the other hand, in relation to our
Corollaries\rl{coro1} and\rl{coro2}, we also refer the reader to Theorem 4.4 and Theorem 4.5 in \cite{CL} where
Caminha and de Lima have recently found another interesting characterizations of the totally umbilical surfaces $L_\tau$
in $\mathcal{H}^{3}$.

\section{An isometrically equivalent model. Steady state type spacetimes}

The steady state space can be expressed in an isometrically equivalent way
as the generalized Robertson-Walker spacetime $-\R \times_{e^t} \R^n$.
That is, the real vector space $\R^{n+1}$ endowed with the Lorentzian
metric tensor
\[
\g{}{}=-dt^2+e^{2t}(dx_1^2+...+dx_n^2),
\]
$(t,x_1,...x_n)$ being the canonical coordinates in $\R^{n+1}=\R\times\R^{n}$. To
see it, take $b \in \R_1^{n+2}$ another null vector such that $\g{a}{b}=1$ and let
$\Phi: \mathcal{H}^{n+1} \rightarrow -\R \times_{e^t}\R^n$ be the map given by \beq
\label{isometria}
\Phi(x)=\left(\log(\g{x}{a}),\frac{x-\g{x}{a}b-\g{x}{b}a}{\g{x}{a}}\right). \eeq
Then it can easily be checked that
\[
(d \Phi)_x(v)=\left(\frac{\g{v}{a}}{\g{x}{a}},
\frac{(v-\g{v}{a}b-\g{v}{b}a)\g{x}{a}
-\g{v}{a}(x-\g{x}{a}b-\g{x}{b}a)}{\g{x}{a}^2} \right)
\]
and that $\Phi$ is an isometry between both spaces which conserves time
orientation.

Therefore, a natural extension of the steady state space consists on
considering the wider family of Lorentzian manifolds defined as follows.
Let $M^n$ be a connected $n$-dimensional Riemannian manifold and consider
the product manifold $\R \times M^n$ endowed with the Lorentzian metric
tensor
\[
\g{}{}=-\pi^\ast_{\mathbb{R}}(dt^2)+e^{2t} \pi^\ast_M( \g{}{}_M)
\]
where $\pi_{\mathbb{R}}$ and $\pi_M$ denote the projections from $\R
\times M^n$ onto each factor, and $\g{}{}_M$ is the Riemannian metric on
$M$. For simplicity, we will write
\[
\g{}{}=-dt^2+e^{2t}\g{}{}_M.
\]
We will denote by \Nn\ the $(n+1)$-dimensional product manifold $\R \times
M^n$ endowed with that Lorentzian metric, and we will refer to them as
steady state type spacetimes. For instance, when $M^n$ is the flat $n$-torus we get the de Sitter
cusp as defined in \cite{Ga}.

\section{Spacelike hypersurfaces in a steady state type spacetime}
Let $\psi:\Sigma^n \fle \Nn$ be a spacelike hypersurface. Observe that
\[
\partial_t=(\partial/\partial t)_{(t,x)}, \,\, t \in \R,\, x \in M,
\]
is a unitary timelike vector field globally defined on the ambient
spacetime which determines a time-orientation on \Nn. Then, the
future-pointing Gauss map of $\Sigma$ is the unique unitary timelike
normal vector field $N$ globally defined on $\Sigma$ in the same
time-orientation as $\partial_t$. Thus, we have
\[
\g{N}{\partial_t} \leq -1 < 0 \quad \mathrm{on} \quad \Sigma.
\]
We will denote by $\Theta:\Sigma \fle (-\infty,-1]$ the smooth function on
$\Sigma$ given by $\Theta=\g{N}{\partial_t}$. Observe that the function
$\Theta$ measures the hyperbolic angle $\theta$ between the
future-pointing vector fields $N$ and $\partial_t$ along $\Sigma$. Indeed,
they are related by $\cosh {\theta}=-\Theta$.

For a spacelike hypersurface $\psi:\Sigma^n \fle \Nn$, we define the height
function of $\Sigma$, denoted by $h$, as the projection of $\Sigma$ onto
$\R$, that is, $h \in \fs$ is the smooth function given by
$h=\pi_{\mathbb{R}}\circ\psi$. In the particular case when $M^n=\R^n$,
by \eqref{isometria} the
height function of $\Sigma$ can be written in terms of the function $u$ as
\[
h(p)=\log{u(\Phi^{-1}(\psi(p)))}, \quad p \in \Sigma.
\]
Therefore, the umbilical hypersurfaces $L_{\tau}$, for which the
function $u$ takes the constant value $u=\tau$, correspond to the slices
$\{\mathrm{log} (\tau) \}\times\R^n $. Actually, for a general steady state type spacetime \Nn, each leaf of the foliation
$\{t\}\times M$ (called here a \textit{slice}) of \Nn\ is a totally umbilical spacelike
hypersurface with constant mean curvature $H=1$.

According to the terminology introduced for hypersurfaces in the steady state space, we will say that a spacelike
hypersurface
$\Sigma$ in \Nn\ is bounded away from the future infinity if there exists $\overline{t}\in\R$ such that
\[
\psi(\Sigma)\subset \{(t,x) \in \Nn: t\leq\overline{t}\},
\]
and we will say that $\Sigma$ is bounded away from the past infinity if there exists $\underline{t}\in\R$ such that
\[
\psi(\Sigma)\subset \{(t,x) \in \Nn: t\geq\underline{t}\}.
\]
We will say that $\Sigma$ is bounded away from the infinity if it is bounded away both from the past and the future
infinity. The following result is the natural extension of our Lemma\rl{lemma}.
\begin{lemma}
\label{lemmatopo}
Let $M^n$ be a Riemannian manifold. If \Nn\ admits a
complete spacelike hypersurface $\Sigma$ which is bounded away from the future infinity, then $M$ is necessarily
complete and the projection $\Pi=\pi_M \circ\psi:\Sigma\fle M$ is a covering map.
\end{lemma}
\begin{proof}
The proof follows the ideas of \cite[Lemma 3.1]{ARS1}. Let $\psi:\Sigma^n\fle \Nn$ be a complete spacelike hypersurface
and let $\Pi=\pi_M \circ \psi:\Sigma \fle M$ denote its projection on $M$. Observe that
\[
\g{v}{v}=-\g{d\psi_p(v)}{\partial_t}^2+e^{2h(p)}\g{d\Pi_p(v)}{d\Pi_p(v)}_M\leq e^{2h(p)}\g{d\Pi_p(v)}{d\Pi_p(v)}_M
\]
for every $p\in\Sigma$ and $v\in T_p\Sigma$. Thus, $\Pi^\ast(\g{}{}_M) \geq (1/e^{2h}) \g{}{}$.
As $\Sigma$ is bounded away from the future infinity, it follows from here that
\[
\Pi^\ast(\g{}{}_M) \geq\frac{1}{e^{2\overline{t}}}\g{}{}=\widetilde{\g{}{}}
\]
for a certain real number $\overline{t}$. Then
reasoning as in the
proof of Lemma\rl{lemma}, we get that $\Pi:(\Sigma,\widetilde{\g{}{}})\fle(M,\g{}{}_M)$ is a local diffeomorphism which
increases the
distance. The proof finishes recalling that if a map, from a connected complete Riemannian manifold
$M_1$ into another connected Riemannian manifold $M_2$ of the same
dimension, increases the distance, then it is a covering map and $M_2$ is
complete \cite[Chapter VIII, Lemma 8.1]{KN}.
\end{proof}

Let $\psi:\Sigma^n \fle \Nn$ be a spacelike hypersurface with height function $h$.
Observe that the gradient of $\pi_{\mathbb{R}}$ on \Nn\ is
\[
\nablabar{\pi_{\mathbb{R}}}=-\g{\nablabar{\pi_{\mathbb{R}}}}{\partial_t}\partial_t=-\partial_t.
\]
Therefore, the gradient of $h$ on $\Sigma$ is
\[
\nabla h=(\nablabar{\pi_{\mathbb{R}}})^\top=-\partial_t^\top,
\]
where $\partial_t^\top$ denotes the tangential component of $\partial_t$
along $\Sigma$, that is,
\[
\partial_t=\partial_t^\top-\Theta N.
\]
In particular,
\beq
\label{normgradh}
\| \nabla h \|^2=\Theta^2-1.
\eeq

Let $Y$ denote the timelike vector field on \Nn\ given by
\[
Y(t,x)=e^t(\partial/\partial t)_{(t,x)}, \qquad (t,x)\in\Nn.
\]
From the relation between $\nablabar$, the Levi-Civita connection of the ambient space \Nn, and
the Levi-Civita connection of $M$, we get that
\beq
\label{eqaux}
\nablabar_ZY=e^t Z
\eeq
for any vector field $Z$ on \Nn. That is, $Y$ is a globally defined closed conformal vector field on \Nn. Writing
$\partial_t=-\nabla h-\Theta N$ along the hypersurface $\Sigma$ and using
Gauss  and Weingarten formulae, it is not difficult to get from
\eqref{eqaux} that
\[
e^h X= \g{\nabla e^h}{X} \partial_t - e^h \nabla_X \nabla
h+e^h\g{AX}{\nabla h}N+e^h \Theta AX-e^h X(\Theta) N
\]
for every $X \in T\Sigma$. Therefore,
\[
\nabla_X\nabla h=\Theta AX-X-\g{\nabla h}{X}\nabla h
\]
for every $X \in T\Sigma$, and the Laplacian on $\Sigma$ of the height function is given by
\beq
\label{laplah}
\Delta h=-n H \Theta -(n+\|\nabla h\|^2).
\eeq

Now we are ready to state the main result of this section, which is the natural extension of our Theorem \ref{deSitter}
to the wider family of steady state type spacetimes.
\begin{theorem}
\label{thgen}
Let $M^n$ be a (necessarily complete) Riemannian manifold
with non-negative sectional curvature and let $\psi:\Sigma^n\fle\Nn$ be a
complete spacelike hypersurface with constant mean curvature $H$. If $\Sigma$ is bounded away from the infinity
then $H=1$. Moreover, when $n=2$, the surface $\Sigma$ is necessarily a slice $\{t\} \times M$.
\end{theorem}
\begin{proof}
The Gauss equation of $\Sigma$ is given by
\[
\g{R(X,Y)X}{Y}=\g{\overline{R}(X,Y)X}{Y}-\g{AX}{X}\g{AY}{Y}+\g{AX}{Y}^2
\]
for every $X,Y \in T\Sigma$, where $R$ and $\overline{R}$ stand for the curvature tensor of $\Sigma$ and $\Nn$,
respectively. Taking traces here, we get that for every $X \in T\Sigma$, $\|X\|=1$,
\begin{eqnarray}
\label{RRM1}
\Ric(X,X) & = & \sum_{i=1}^m\g{\overline{R}(X,E_i)X}{E_i}+nH\g{AX}{X}+\g{AX}{AX}\\
\nonumber {} & \geq & \sum_{i=1}^m\g{\overline{R}(X,E_i)X}{E_i}-\frac{n^2H^2}{4},
\end{eqnarray}
where $\Ric$ stands for the Ricci curvature of $\Sigma$. On the other hand, taking into account the relationship between
the curvature tensor of $\Nn$ and the curvature tensor $R_M$ of $M$ (as well as the derivatives of the warping function)
\cite[Proposition 42]{ONe}, we also have that
\beq
\label{RRM}
\g{\overline{R}(X,E_i)X}{E_i}=e^{2h}\g{R_M(X^\ast,E_i^\ast)X^\ast}{E_i^\ast}_M-\g{X}{E_i}^2+1,
\eeq
for $1 \leq i \leq n$, where $X^\ast=(\pi_M)_\ast(X)$ denotes the projection of
the vector field $X \in T\Sigma$ onto $M$, that is, $X=X^\ast-\g{X}{\partial_t}\partial_t$.
Denoting by $K_M(X^\ast\wedge E_i^\ast)$ the sectional curvature in $M$ of the 2-plane generated by
$X^\ast$ and $E_i^\ast$, we have that
\[
\g{R_M(X^\ast,E_i^\ast)X^\ast}{E_i^\ast}_M=K_M(X^\ast\wedge E_i^\ast)\|X^\ast\wedge E_i^\ast\|_M^2,
\]
where, as usual,
\[
\|X^\ast\wedge E_i^\ast\|_M^2=\g{X^\ast}{X^\ast}_M \g{E_i^\ast}{E_i^\ast}_M - \g{X^\ast}{E_i^\ast}_M^2.
\]
Using this into \rf{RRM}, we get from \rf{RRM1}
\begin{eqnarray}
\label{Ric2}
\Ric(X,X) & \geq & e^{2h}\sum_{i=1}^nK_M(X^\ast\wedge E_i^\ast)\|X^\ast\wedge E_i^\ast\|_M^2+n-1-\frac{n^2H^2}{4}\\
\nonumber {} & \geq & n-1-\frac{n^2H^2}{4}.
\end{eqnarray}

Consequently, the assumption on the sectional curvature of $M$ implies
that the Ricci curvature tensor of $\Sigma$ is bounded from below. The
proof that $H=1$ follows now by applying Lemma\rl{OmoriYau} to the height function.
Actually, since $h$ is bounded from above we know that there exists a sequence
$\{p_k\}_{k\in\mathbb{N}}$ in $\Sigma$ such that
\[
\lim_{k \rightarrow \infty} h(p_k)=\sup_{\Sigma} h<+\infty, \quad
\|\nabla h (p_k)\|^2=\Theta^2(p_k)-1<1/k^2
\]
and
\[
\Delta h(p_k)=-n H \Theta(p_k)-(n+\|\nabla h(p_k)\|^2)<1/k.
\]
Therefore $\lim_{k \rightarrow \infty} \Theta(p_k)=-1$, and taking limits in the last
equation we get $H \leq 1$.

In a similar way, by applying Lemma \ref{OmoriYau} to the function $-h$ we obtain another sequence $\{q_k\}_{k\in\mathbb{N}}$ in
$\Sigma$ such that
\[
\lim_{k \rightarrow \infty} h(q_k)=\inf_{\Sigma} h>-\infty, \quad
\|\nabla h (q_k)\|^2=\Theta^2(q_k)-1<1/k^2
\]
and
\[
\Delta h(q_k)=-n H \Theta(q_k)-(n+\|\nabla h(q_k)\|^2)>-1/k.
\]
Thus $\lim_{k \rightarrow \infty} \Theta(q_k)=-1$ again, and taking limits in the last
equation we get $H\geq 1$. As a consequence, $H=1$.

Consider now the $2$-dimensional case. Since $H=1$, \rf{Ric2} implies that the Gaussian curvature of $\Sigma$ is
non-negative and, therefore, $\Sigma$ is a parabolic surface. Moreover, \rf{laplah}, jointly with \rf{normgradh}, implies
that
\[
\Delta h=-2 \Theta - (2+\|\nabla h\|^2)=-2 \Theta - (1+ \Theta^2)=-(\Theta+1)^2 \leq 0.
\]
Therefore, $h$ is a superharmonic function bounded from below on the parabolic surface $\Sigma$ and hence it must be
constant.
\end{proof}

\bibliographystyle{amsplain}

\end{document}